\documentclass[11pt, a4paper]{amsart}

\usepackage{amsmath}
\usepackage{amstext}
\usepackage{amssymb}
\usepackage{amsthm}
\usepackage{enumerate}
\usepackage{enumitem}
\usepackage{mathrsfs}
\usepackage[all]{xy}
\usepackage{datetime}
\usepackage[svgnames]{xcolor}
\usepackage{tikz}

\pgfdeclarelayer{edgelayer}
\pgfdeclarelayer{nodelayer}
\pgfsetlayers{edgelayer,nodelayer,main}

\tikzstyle{none}=[inner sep=0pt]

\makeatletter
\def\imod#1{\allowbreak\mkern10mu({\operator@font mod}\,\,#1)}
\makeatother

\newtheorem{theorem}{Theorem}[section]

\newtheorem{prop}[theorem]{Proposition}

\newtheorem{corollary}[theorem]{Corollary}

\theoremstyle{definition}
\newtheorem{definition}[theorem]{Definition}

\theoremstyle{remark}

\theoremstyle{remark}

\numberwithin{equation}{section}

    \DeclareMathOperator{\dom}{dom}
    
    \DeclareMathOperator{\ran}{ran}

    \DeclareMathOperator{\fin}{FIN}

    \newcommand{\restrict}{\!\upharpoonright\!}

    \newcommand{\forces}{\Vdash}

\def\R{{\mathbb R}}

\def\N{{\mathbb N}}

\title{The Ramsey property implies no mad families}

\author{David Schrittesser}

\address{Kurt G\"odel Research Center, University of Vienna, W\"ahringer Strasse 25, 1090 Vienna, Austria}

\email{asgert@math.ku.dk}

\author{Asger T\"ornquist}

\address{Department of Mathematical Sciences, University of Copenhagen, Universitetsparken 5, 2100 Copenhagen, Denmark}

\email{david@logic.univie.ac.at}



\date{\today}


\begin{document}

\begin{abstract}
We show that if all collections of infinite subsets of $\N$ have the Ramsey property, then there are no infinite maximal almost disjoint (mad) families. This solves a long-standing problem going back to Mathias \cite{mathias}. The proof exploits an idea which has its natural roots in ergodic theory, topological dynamics, and invariant descriptive set theory: We use that a certain function associated to a purported mad family is invariant under the equivalence relation $E_0$, and thus is constant on a ``large'' set. 
Furthermore we announce a number of additional results about mad families relative to more complicated Borel ideals.
\end{abstract}

\maketitle

\section{Introduction}

In his seminal paper \cite{mathias}, Mathias established a connection between three different ideas in mathematics: The combinatorial set theory of \emph{maximal almost disjoint families}, infinite dimensional Ramsey theory, and Cohen's method of forcing. He asked if the combinatorial statement ``all sets have the Ramsey Property'' implies that there are no infinite mad families. In this paper we answer this in the affirmative, working in the theory ZF+DC+R-Unif (See Definition \ref{def2}).

Let us recall the key notions: An \emph{almost disjoint} family (on the natural numbers $\N$) is a family $\mathcal A$ of \emph{infinite} subsets of $\N$ such that if $x,y\in\mathcal A$ then either $x=y$, or $x\cap y$ is finite. A \emph{maximal almost disjoint family} (``mad family'') is an almost disjoint family which is not a proper subset of an almost disjoint family. Finite mad families are easily seen to exists, e.g. $\{E,O\}$, where $E$ is the set of even numbers and $O$ is the set of odd numbers. The existence of \emph{infinite} mad families follows easily from Zorn's lemma (equivalently, the Axiom of Choice).

Given a set $X$ and a natural number $k\in\N$, let $[X]^k$ denote the set of all subsets of $X$ with exactly $k$ elements. The classical infinite Ramsey theorem in combinatorics says if $S\subseteq [\N]^k$ then there is an infinite set $B\subseteq \N$ such that either $[B]^k\subseteq S$, or $[B]^k\cap S=\emptyset$. Motivated by a question of Erd\H{o}s and Rado, infinite dimensional generalizations of this theorem were discovered in the 1960s and 1970s. In this paper, we denote by $[X]^\infty$ the set of \emph{countably infinite} subsets of $X$. Moreover, given $S\subseteq [\N]^\infty$, we will say that $S$ has the \emph{Ramsey Property}, or simply \emph{is Ramsey}, if there is $B\in [\N]^\infty$ such that either $[B]^\infty\subseteq S$, or $[B]^\infty\cap S=\emptyset$. Erd\H{o}s  and Rado showed that the Axiom of Choice implies that not all sets $S\subseteq [\N]^\infty$ are Ramsey. Later, in \cite{galvin-prikry} and \cite{silver}, it was shown that Borel and analytic $S\subseteq [\N]^\infty$ \emph{are} Ramsey, and finally Ellentuck in 1973 \cite{ellentuck} characterized the Ramsey Property in terms of Baire measurability in the \emph{Ellentuck topology} on $[\N]^\infty$.

Concurrent with these developments in Ramsey theory, Cohen's introduction of the method of forcing for independence proofs in set theory in the early 1960s set off an explosion of independence results, among the most famous of which is Solovay's model of Zermelo-Fraenkel set theory (ZF) in which only a fragment of the Axiom of Choice--namely Dependent Choice (DC)---holds and in which all subsets of $\R$ are Lebesgue and have the Baire Property. It was Mathias, in his \emph{Happy Families} paper (drafts of which circulated already in the late 1960s), who connected forcing to Erd\H os and Rado's question, and to mad families. Mathias did this by introducing what is now known as \emph{Mathias forcing}, which he used to show that in Solovay's model \emph{all} sets $S\subseteq [\N]^\infty$ are Ramsey. Mad families, and their connection to Mathias forcing, plays a central role in this proof.

Mathias asked two central questions, which his methods did not allow him to answer at the time: (1) Are there infinite mad families in Solovay's model? (2) If all sets $S\subseteq [\N]^\infty$ are Ramsey, does it follow that there are no infinite mad families? A positive answer to (2) would give a negative answer to (1). 

There was only modest progress on these questions until very recently, when suddenly the research in mad families and forcing experienced a renaissance. Question (1) was solved in 2014 in \cite{tornquist}, and shortly after, Horowitz and Shelah showed in \cite{horowitz-shelah-can} that a model of ZF in which there are no mad families can be achieved without using an inaccessible cardinal, which is otherwise a crucial ingredient in the construction of Solovay's model. Neeman and Norwood in \cite{neeman-norwood} and independently, Bakke Haga in joint work with the present authors in \cite{haga-schrittesser-tornquist} proved a number of further results, among them that $V=L(\R)$+AD implies there are no mad families. Horowitz and Shelah also solved a number of related questions that had been formulated over the years, in particular, they showed the existence of a Borel ``med'' family in \cite{horowitz-shelah}, see also \cite{schrittesser} for a simpler proof.

We denote by R-Unif the principle of uniformization on Ramsey positive sets (see see Definition \ref{def2} below; Solovay's model easily satisfies this principle). In this paper we give the following positive solution to Mathias' question.

\begin{theorem}\label{t.main1}
(ZF+DC+R-Unif) If all sets have the Ramsey Property then there are no infinite mad families.
\end{theorem}

We note that this implies the main results of \cite{tornquist} and \cite{neeman-norwood}.

Theorem \ref{t.main1} may seem all the more surprising given another recent result of Horowitz and Shelah \cite{horowitz-shelah-madness}, who show that for a variety of measurability notions including Lebesgue measure, ``all sets are  measurable'' is compatible with the existence of an infinite mad family.

\medskip

Let us briefly comment on the proof of Theorem \ref{t.main1}, and the difficulties that have to be overcome. For this discussion, suppose $\mathcal A\subseteq\mathcal [\N]^\infty$ is an infinite mad family, and assume ``all sets are Ramsey'' and ``Ramsey uniformization'' (again see Definition \ref{def2}). The first difficulty encountered is that the set of $x\in [\N]^\infty$ which meet exactly one element of $\mathcal A$ in an infinite set is clearly Ramsey co-null when $\mathcal A$ is a mad family. We will overcome this by associating to each $z\in [\N]^\infty$ a very, very sparse set $\tilde z\in [\N]^\infty$, which is constructed using a fixed, infinite sequence $(x_n)_{n\in\N}$ chosen from $\mathcal A$ (it is here that we use the Principle of Dependent Choice, i.e., DC). A basic property of the map $z\mapsto \tilde z$ is that it is equivariant under finite differences, that is, if $z\triangle z'$ is finite then $\tilde z\triangle\tilde z'$ is finite.

Because we assumed that $\mathcal A$ is maximal, for each $z\in [\N]^\infty$ there is some $y_z\in\mathcal A$ such that $\tilde z\cap y$ is infinite, and so R-Unif gives us a function $f:[\N]^\infty\to\mathcal A$ such that $f(z)\cap \tilde z$ is infinite for $z$ in a Ramsey positive set. The special way that $z\mapsto \tilde z$ will be defined below will ensure that no uniformizing function $f$ can have the invariance property that if $|z\triangle z'|<\infty$ implies $f(z)=f(z')$. While there is no reason to expect that an abstract application of R-Unif would give us $f$ with this property, it turns out that with some work we can get dangerously close to having such an invariant $f$. Indeed, by using the assumption that all sets are Ramsey we can find an infinite set $W\subseteq\N$ such that the restriction $f\restrict [W]^\infty$ is continuous, and so the range $f([W]^\infty)$ is an analytic set. Using that $f([W]^\infty)$ is analytic, we will define a function $z\mapsto T^z$, where $T^z$ can be thought of as a tree of approximations to possible, natural uniformization functions. It then turns out that the map $z\mapsto T^z$ satisfies that if $|z\triangle z'|<\infty$ then $T^z=T^{z'}$. This in turn leads to that $z\mapsto T^z$ is constant on a Ramsey positive set, which then leads to a contradiction.

\subsection*{Acknowledgements} The first author thanks the Austrian Science Fund FWF for support through grant P29999. The second author thanks the Danish Council for Independent Research for generous support through grant no.\ 7014-00145B. The second author also gratefully acknowledges his association with the Centre for Symmetry and Deformation, funded by the Danish National Research Foundation (DNRF92).

\section{Notation and background definitions}

In this section we summarize the background needed for the proof. A good general reference for all the background needed is \cite{kechris}. A comprehensive treatise on modern, infinitary Ramsey theory can be found in \cite{todorcevic}.

\subsection{Descriptive set theory}

A topological space $X$ is called \emph{Polish} if it is separable and admits complete metric that induces the topology. In this paper we will be working with the Polish space $2^\N=\{0,1\}^\N$ and $\N^\N$ (with the product topology, taking $\{0,1\}$ and $\N$ discrete) and subspaces of these space. Recall the following key notion from descriptive set theory:

\begin{definition}
A subset $A\subseteq X$ of a Polish space $X$ is \emph{analytic} if there is a continuous $f:Y\to X$ from a Polish space $Y$ to $X$ such that $A=\ran(f)$.
\end{definition}

Since $\N^\N$ maps continuously onto any Polish space we have: \emph{$A\subseteq X$ is analytic iff there is a continuous $f:\N^\N\to X$ such that $\ran(f)=A$}. We will use this characterization as our definition of analytic set below.

For the proof of Theorem \ref{t.main1} we need the following combinatorial description of the topology on $\N^\N$. We denote by $\N^n$ the set of all functions $s:\{1,\ldots, n\}\to\N$, and we let $\N^{<\N}=\{\emptyset\}\cup\bigcup_{n\in\N} \N^n$. (We shall think of $\emptyset$ as the function with empty domain, which is why it is included as an element of $\N^{<\N}$). For $s,t\in\N^{<\N}\cup\N^\N$ we will write $s\subseteq t$ (``t extends s'') if $\dom(s)\subseteq\dom(t)$ and $s(i)=t(i)$ for all $i\in\dom(s)$; we will write $s\perp t$ (``$s$ and $t$ are incompatible'') if $s\not\subseteq t$ and $t\not\subseteq s$.

For each $s\in\N^{<\N}$, let
$$
N_s=\{f\in\N^\N: (\forall i\in\dom(s))\ f(i)=s(i)\}.
$$
The family $\{N_s: s\in\N^{<\N}\}$ is easily seen to form a basis for the topology on $\N^\N$.

Note that $\N^{<\N}$ is countable, and so $2^{\N^{<\N}}$ is a Polish space (isomorphic to $2^\N$) in the product topology, taking $2=\{0,1\}$ discrete. This view will be important later in the proof of Theorem \ref{t.main1} where we will describing the properties of a certain continuous function $f$ defined on $\N^\N$ in terms of a ``derived'' function $f':\N^\N\to 2^{\N^{<\N}}$.

\subsection{The Ramsey Property} For any set $X$ we define
$$
[X]^{\infty}=\{A\subseteq X: A\text{ is infinite}\}.
$$
Recall from the introduction that a set $S\subseteq\ [\N]^\infty$ is \emph{Ramsey} (or \emph{has the Ramsey property} if there is $B\in [\N]^\infty$ such that $[B]^\infty\subseteq S$ or $S\cap [B]^\infty=\emptyset$.

To each infinite $B\subseteq\N$ we let $\hat B:\N\to\N$ be the unique increasing function with $B=\ran(\hat B)$. For each infinite $A\subseteq\N$ the map $[A]^\infty\to\N^\N: B\mapsto \hat B$ naturally identifies $[A]^\infty$ with a subset of $\N^\N$. The (subspace) topology that $[A]^\infty$ inherits under this identification will be called the \emph{standard topology} on $[A]^\infty$.

The crucial use in our proof of the Ramsey Property comes from the following well-known fact:

\begin{prop}\label{p.continuous}
(ZF+DC) Suppose all sets have the Ramsey Property and let $\vartheta:[\N]^{\N}\to Y$ be a function, where $Y$ is a separable topological space.  Then there is $A\in [\N]^\infty$ such that $\vartheta\restrict [A]^{\infty}$ is continuous with respect to the standard topology on $[A]^\infty$.
\end{prop}

\subsection{Uniformization} We need the following \emph{uniformization} principle.

\begin{definition}\label{def2}
(1) Let $X,Y$ be Polish spaces and let $R\subseteq X\times Y$ be set with the property that $X=\dom(R)$, where
$$
\dom(R)\overset{\text{def}}=\{x\in X: (\exists y\in Y)\ (x,y)\in R\}.
$$ 
If $Z\subseteq X$ and $\vartheta: Z\to Y$, we say that $\vartheta$ \emph{uniformizes} $R$ on $Z$ if for all $x\in Z$ we have $(x,\vartheta(x))\in R$.

(2) The \emph{Ramsey uniformization principle}, abbreviated R-Unif, is the following statement: For all $R\subseteq [\N]^{\infty}\times Y$, where $Y$ is Polish, and all infinite $A\subseteq\N$ such that $[A]^\infty\subseteq\dom(R)$ there is $B\in [A]^\infty$ and $\vartheta: [B]^\infty\to Y$ which uniformizes $R$ on $[B]^\infty$.
\end{definition}

That Ramsey uniformization property holds 
in Solovay's model follows by the argument given in \cite[1.12, p. 46]{solovay} with Random forcing replaced everywhere by Mathias forcing.

\subsection{Invariance under finite change} The last ingredient for the proof is the notion of ``$E_0$-invariance''.

\begin{definition}
(1) $E_0$ is the equivalence relation defined on $[\N]^\infty$ by
$$
x E_0 y\iff |x\triangle y|<\infty.
$$
In other words, $x,y\in [\N]^\infty$ are $E_0$ equivalent iff they differ only on a finite set.

(2) A function $f:[\N]^\infty\to Y$ is called \emph{$E_0$-invariant} if $x E_0 y$ implies $f(x) E_0 f(y)$.
\end{definition}

\section{Proof of Theorem 1}

We work under the following assumptions: ZF+\-DC+\-R-Unif+ ``All sets have the Ramsey Property''. 

Let $\mathcal A\subseteq [\N]^\infty$ be an infinite almost disjoint family. We will show that $\mathcal A$ is not maximal. 

Let $(x_n)_{n\in\N}$ be an injective sequence of elements in $\mathcal A$ (here we use that DC implies that all infinite sets are Dedekind infinite). We may assume that $\mathcal A'=\mathcal A\setminus\{x_n: n\in\N\}$ is non-empty, since otherwise an easy diagonalization shows that $\mathcal A$ is not maximal.
Moreover, by possibly replacing $x_n$ by $x_n\setminus (\bigcup_{i<n} x_i)$, let us assume for simplicity that if $n\neq m$ then $x_n\cap x_m=\emptyset$.

Recall that when  $z\in [\N]^\infty$ then $\hat z:\N\to\N$ is the unique increasing function such that $\ran(\hat z)=z$. Using the sequence $(x_n)_{n\in\N}$ fixed above, define for each $z\in [\N]^{\infty}$,
$$
\tilde z=\{\widehat{x_{\hat z(n)}}(\hat z(n+1)):n\in\N\}.
$$
Note that $|\tilde z\cap x_n|\leq 1$ for all $n$, so proving the following claim will prove the theorem:

\medskip
\noindent {\bf Main Claim.} There is $z\in [\N]^\infty$ such that for all $y\in\mathcal A'$, $|\tilde z\cap y|<\infty$.
\medskip

\medskip

Suppose the claim is false. Then by Ramsey Uniformization there is $W\in [\N]^{<\infty}$ and $\vartheta: [W]^{\infty}\to \mathcal A'$ such that
$$
|\vartheta(z)\cap\tilde z|=\infty
$$
for all $z\in [W]^\infty$. By Proposition \ref{p.continuous} we may assume (after possibly replacing $W$ with an infinite subset of $W$) that $\vartheta\restrict [W]^\infty$ is continuous. Then $\mathcal B=\vartheta( [W]^{\infty})$ is an analytic subset of $\mathcal A'$, and we fix a continuous function $f:\N^\N\to [\N]^\infty$ such that $\ran(f)=\mathcal B$.

For $z\in [\N]^\infty$, let
$$
T^z=\{s\in \N^{<\N}: (\exists x\in N_s)\ |\tilde z\cap f(x)|=\infty\}.
$$
By identifying $\mathcal P(\N^{<\N})$ with $2^{\N^{<\N}}$, we will think of $z\mapsto T_z$ as a map $[\N]^\infty\to 2^{\N^{<\N}}$. The reader can easily verify that $T^z$ is a tree in the sense of \cite{kechris}, and that $\emptyset\in T^z$ for all $z\in [W]^\infty$.

\medskip

\noindent {\bf Subclaim 1.} The function $z\mapsto T^z$ is $E_0$-invariant

\medskip

{\it Proof.} Suppose $|z'\triangle z|<\infty$. Then we can find $k_0,k_0'$ such that $\hat z(k_0+i)=\hat z'(k_0'+i)$ for all $i\in\N$. Then 
$$
\widehat{x_{\tilde z(k_0+i)}}(\hat z(k_0+i+1))=\widehat{x_{\tilde z'(k_0'+i)}}(\hat z'(k_0'+i+1))
$$
for all $i\in\N$. It follows that $\tilde z E_0 \tilde z'$, but then $|\tilde z\cap f(x)|=\infty$ if and only if $|\tilde z'\cap f(x)|=\infty$, so $T^z=T^{z'}$.\hfill{\tiny Subclaim 1.}$\dashv$

\medskip

\noindent  {\bf Subclaim 2.} There is $W_0\in [W]^\infty$ such that $z\mapsto T^z$ is constant on $[W_0]^\infty$.

\medskip

{\it Proof.} By proposition \ref{p.continuous} there is $W_0\in [W]^\infty$ such that $z\mapsto T^z$ is continuous on $[W_0]^\infty$. Since $z\mapsto T^z$ is $E_0$ invariant, it follows that $z\mapsto T^z$ is constant on $[W_0]^\infty$.\hfill {\tiny Subclaim 2.}$\dashv$

\medskip

From now on we fix $W_0\in [W]^\infty$ and $\tilde T\subseteq \N^{<\N}$ such that $T^z=\tilde T$ for all $z\in [W_0]^\infty$. The next claim echoes the claim on top of p. 65 in \cite{tornquist}.

\medskip

\noindent  {\bf Subclaim 3.} Suppose there are $t^0,t^1\in\
\tilde T$ and $n_0\in\N$ such that for all $y_0\in f(N_{t^0})$ and $y_1\in f(N_{t^1})$ we have $n_0\in f(y_0)\triangle f(y_1)$. Then there are $s^0,s^1\in \tilde T$ and $k\in\N$ such that $s^0\supseteq t^0$, $s^1\supseteq t^1$, and for all $y_0\in f(N_{s_0})$ and $y_1\in f(N_{s_1})$ we have $y_0\cap y_1\subseteq \{1,\ldots, k\}$.

\medskip

{\it Proof.} Suppose \emph{no} $s^0\supseteq t^0$ and $s^1\supseteq t^1$, with $s^0,s^1\in\tilde T$, satisfies the claim. Then for every $m\in\N$, $t^0\subseteq u\in \tilde T$ and $t^1\subseteq v\in\tilde T$ we can find $m'>m$, $u\subseteq u'\in\tilde T$ and $v\subseteq v'\in\tilde T$ such that for some $x_0\in N_{u'}$ and $x_1\in N_{v'}$ we have $m'\in f(x_0)\cap f(x_1)$. By the continuity of $f$ we can then find $u'\subseteq u''\in\tilde T$ and $v'\subseteq v''\in \tilde T$ such that \emph{for all} $x_0\in N_{u''}$ and $x_1\in N_{v''}$ we have $m'\in f(x_0)\cap f(x_1)$.

Using the previous paragraph repeatedly, we can now build sequences 
\begin{align*}
&t^0\subseteq u_1\subseteq u_2\subseteq\cdots\\
&t^1\subseteq v_1\subseteq v_2\subseteq\cdots\\
&m_0<m_1<\cdots
\end{align*}
where $u_i,v_i\in \tilde T$, and for all $x_0\in N_{u_i}$ and $x_1\in N_{v_i}$ we have $m_i\in f(x_0)\cap f(x_1)$ when $i>0$. Let $u_\infty, v_\infty\in\N^\N$ be such that $u_i\subseteq u_\infty$ and $v_i\subseteq v_\infty$ for all $i\in\N$. Then $|f(u_\infty)\cap f(v_\infty)|=\infty$ since $m_i\in f(u_\infty)\cap f(v_\infty)$ for all $i> 0$, but $f(u_\infty)\neq f(v_\infty)$ since $n_0\in f(u_\infty)\triangle f(v_\infty)$. This contradicts that $\ran(f)$ is a subset of the almost disjoint family $\mathcal A$. \hfill {\tiny Subclaim 3.}$\dashv$

\medskip

\noindent  {\bf Subclaim 4.} There is a unique $y^*\in\mathcal B$ such that for all $z\in [W_0]^\infty$ we have $|\tilde z\cap y^*|=\infty$.

\medskip

{\it Proof:} Since  $\vartheta(z)\in\mathcal B$ for $z\in [W_0]^\infty$ and $|\tilde z\cap \vartheta(z)|=\infty$ by definition, for every $z\in [W_0]^\infty$ there is some $y\in\mathcal B$ such that $|\tilde z\cap y|=\infty$. We must show that there is a unique $y\in\mathcal B$ not depending on $z$ satisfying this.

Suppose not, and let $x_0,x_1\in\N^\N$ such that $f(x_0)\neq f(x_1)$ and for some $z_0,z_1\in [W_0]^\infty$ we have $|\tilde z_0\cap f(x_0)|=\infty$ and $|\tilde z_1\cap f(x_1)|=\infty$. By continuity of $f$ we can find $t_0\subseteq x_0$ and $t_1\subseteq x_1$ and $n_0\in\N$ such that for all $x\in N_{t^0}$ and $x'\in N_{t^1}$ we have $n_0\in f(x)\triangle f(x')$. 

Note that $t^0,t^1\in \tilde T$ since $t^i\in T^{z_i}$ and $T^{z_i}=\tilde T$ for $i\in\{0,1\}$. By the previous subclaim we can find $t^0\subseteq s^0\in\tilde T$ and $t^1\subseteq s^1\in\tilde T$ and $k\in\N$ such that for all $y_0\in f(N_{s^0})$ and $y_1\in f(N_{s^1})$ we have $y_0\cap y_1\subseteq \{1,\ldots, k\}$.

Let
$$
x_n^i=x_n\cap \bigcup f(N_{s^i}).
$$
By our assumptions on $s^0$ and $s^1$ we have that $x_n^0\cap x_n^1=\emptyset$ for $n$ sufficiently large. By possibly removing a finite initial segment from $W_0$, we may assume that $x_n^0\cap x_n^1=\emptyset$ \emph{for all} $n\in W_0$.

\medskip

Below, for $A\subseteq\N$, we let $A\setminus n=\{i\in A: i>n\}$. Clearly, for each $n\in W_0$ at least one of the following hold:
\begin{enumerate}[label=(\arabic*)]
\setcounter{enumi}{-1}
\item $(\exists^\infty j\in W_0\setminus n)\  \hat x_n(j)\in x_n^0$;
\item $(\exists^\infty j\in W_0\setminus n)\  \hat x_n(j)\notin x_n^0$.
\end{enumerate}

By refining $W_0$ to $W_{0}'\in [W_0]^\infty$ we can then arrange (keeping in mind that $x_n^0\cap x_n^1=\emptyset$) that for each $n\in W_{0}'$ exactly one of the following hold:
\begin{enumerate}[label=(\arabic*$'$)]
\setcounter{enumi}{-1}
\item $(\forall j\in W_{0}'\setminus n)\  \hat x_n(j)\in x_n^0$;
\item $(\forall j\in W_{0}'\setminus n)\  \hat x_n(j)\notin x_n^0$.
\end{enumerate}

By refining $W_0'$ one more time to $W_{0}''\in [W_0']^\infty$ we can then arrange that exactly one of the following hold:

\begin{enumerate}[label=(\arabic*$''$)]
\setcounter{enumi}{-1}
\item $(\forall n\in W_0'')(\forall j\in W_{0}''\setminus n)\  \hat x_n(j)\in x_n^0$;
\item $(\forall n\in W_0'')(\forall j\in W_{0}''\setminus n)\  \hat x_n(j)\notin x_n^0$.
\end{enumerate}

Now we arrive at a contradiction: If the (0$''$) holds, then since $x_n^0\cap x_n^1=\emptyset$ for all $n\in W_0''$ we have for all $z\in [W_0'']^\infty$ that $\tilde z\cap\bigcup f(N_{s^1})=\emptyset$, contradicting that $s^1\in\tilde T=T^z$. Similarly, if (1'') holds we get for all $z\in [W_0'']^\infty$ that $\tilde z\cap \bigcup f(N_{s^0})=\emptyset$, contradicting that $s^0\in\tilde T$.\hfill {\tiny Subclaim 4.}$\dashv$

\medskip

To finish the proof of the Main Claim, let $y^*\in\mathcal B$ be as in the previous claim. Since $y^*\in\mathcal A'$ we have $x_n\cap y^*$ is finite for all $n\in\N$. Let $z\in [W_0]^\infty$ be such that
$$
\hat x_{\hat z(n)}(\hat z(n+1))>\max(x_{\hat z(n)}\cap y^*)
$$
for all $n\in\N$. Then $\tilde z\cap y=\emptyset$, contradicting Subclaim 4. This contradiction establishes the Main Claim, and as noted above, the Main Claim easily implies that $\mathcal A$ is not maximal, which is what we needed to prove.\hfill $\square$

\medskip

In the proof above, a crucial point was obtaining $W''_0 \in [W_0]^\infty$ such that $(\forall z \in [W''_0]^\infty)\;\tilde z \subseteq x^0_n$ or $\tilde z \cap x^0_n =\emptyset$. 
Note that an alternative and quick way to obtain such $W''_0$ is to appeal to Ramsey's Theorem for pairs and take $W''_0$ to be a homogeneous set for the 2-coloring
\[
c(n,j)=\begin{cases} 1 &\text{ if $\hat x_n(j) \in x^0_n$,}\\
0 &\text{ otherwise.}
\end{cases}
\]

\section{Corollaries and further results}

\begin{corollary}[T\"ornquist \cite{tornquist}]
There are no mad families in Solovay's model.
\end{corollary}

{\it Proof}: By \cite{solovay}, Solovay's model is a model of ZF+DC. That the Ramsey Property holds in this model follows from \cite{mathias}. Finally, that the Ramsey uniformization principle holds by our remarks after Definition \ref{def2}.\hfill $\square$

\medskip

We point out that the proof of Theorem \ref{t.main1} above localizes as follows.
\begin{corollary}
If $\Gamma$, $\Gamma'$ are reasonable pointclasses such that every relation $R\in\Gamma$ can be uniformized on a Ramsey positive set by a $\Gamma'$-measurable function and all sets in $\Gamma'$  are Ramsey, then there are no infinite mad families in $\Gamma$.
\end{corollary}

The local version allows us to draw consequences regarding the Axiom of Projective Determinacy (short PD).

\begin{corollary}[\cite{neeman-norwood}, \cite{haga-schrittesser-tornquist}]
PD implies there are no projective infinite mad families.
\end{corollary}

{\it Proof}:
The hypotheses of the previous theorem hold with $\Gamma'=\Gamma$ equal to the class of projective sets: PD implies that this pointclass has the uniformization property, and by \cite{harrington-kechris}, all projective sets are completely Ramsey under PD.\hfill$\square$

\medskip

Another consequence of our proof is that Mathias forcing destroys mad families from the ground model:
\begin{theorem}
In the Mathias extension, there is no mad family $\mathcal A$ such that 
\begin{itemize}
\item $\mathcal A$ is definable by a $\Sigma_1$ formula in the language of set theory with parameters in the ground model;
\item There is an infinite sequence in the ground model whose elements belong to $\mathcal A$.
\end{itemize}
In particular no infinite a.d.\ family from the ground model is maximal in the Mathias extension.
\end{theorem}
{\it Sketch of proof}:
Work in $V[x]$ where $x$ is Mathias over $V$, and suppose $\mathcal A=\{z \in [\N]^\infty : \Psi(z)\}$ is an infinite a.d.\ family, where $\Psi(z)$ is $\Sigma_1$ (with parameter in $V$).
We show that $\tilde x$ is almost disjoint from every $z \in \mathcal A$, where $\tilde x$ is defined in $V[x]$ as in the proof of Theorem \ref{t.main1} from an infinite sequence in the ground model whose elements belong to $\mathcal A$.
Towards a contradiction, suppose $\tilde x$ is \emph{not} almost disjoint from every $z \in \mathcal A$. Fix a Mathias condition $(s,A)\in V$ with $s \subseteq x \subseteq A$ and a name $\dot y\in V$ such that $p\forces \dot y \cap \tilde{\dot x}_G$ is infinite and $\Psi(\dot y)$ (where ${\dot x}_G$ is a name for the Mathias real).
By a well-known property of Mathias forcing (so-called \emph{continuous reading of names})
we can assume that there is a continuous function $\vartheta\colon[\N]^\infty\to[\N]^\infty$ with code in $V$ such that
$p\forces\dot y= \vartheta(\dot {x}_G)$.
It is easy to see that $\vartheta(y) \in \mathcal A$ for any $y\in V[x]$ such that  $s\subseteq y \in [x]^\infty$ (since such $y$ is also Mathias over $V$; here we also use the definability of $\mathcal A$).
But then in $V[x]$, $\ran(\vartheta)$ would be an analytic almost disjoint family such that any element of
$\{\tilde y : y \in[x]^\infty\}$ has infinite intersection with some element of $\ran(\vartheta)$, which is impossible by the proof of Theorem \ref{t.main1}.\hfill$\square$

\medskip

Surprisingly, the connection with the Ramsey property extends beyond the ideal of finite sets to much more complicated ideals.
We construct such ideals using the familiar Fubini sum:
Given, for each $n\in\N$, an ideal $\mathcal J_n$ on a countable set $S_n$ we obtain an ideal $\mathcal J$ on $S=\bigsqcup_n S_n$ as follows:
\[
\mathcal J = \bigoplus_n \mathcal J_n =\{X\subseteq S : (\forall^\infty n)\; X\cap S_n\in\mathcal J_n\}
\]
where $(\forall^\infty n)$ means ``for all but finitely many $n$.''
The Fubini sum $\bigoplus_n\fin$ (where $\fin$ denotes the ideal of finite sets on $\N$) is also known as $\fin\times\fin$ or $\fin^2$;
iterating Fubini sums into the transfinite we obtain $\fin^ \alpha$, $\alpha < \omega_1$. 
This family of ideals of lies cofinally in the Borel hierarchy in terms of complexity.

The notion of mad family can be extended to arbitrary ideals on a countable set: If $\mathcal J$ is such an ideal, a $\mathcal J$-almost disjoint family is a subfamily $\mathcal A$ of $\operatorname{\mathcal P}(S)\setminus \mathcal J$ such that for any two distinct $A,A'\in \mathcal A$, $A\cap A'\in \mathcal J$. A $\mathcal J$-mad family is of course a $\mathcal J$-almost disjoint family which is maximal under $\subseteq$ among such families.

In the forthcoming article \cite{schrittesser-tornquist} we show the following:
\begin{theorem}(ZF+DC+R-Unif) 
Let $\alpha<\omega_1$. If all sets have the Ramsey Property then there are no infinite $\fin^\alpha$-mad families.
\end{theorem}
As for classical mad families, we immediately obtain corollaries regarding the Axiom of Projective Determinacy, and Solovay"s model.
The first corollary was already shown in \cite{haga-schrittesser-tornquist} using forcing over inner models.
\begin{corollary}(ZF+PD)
For each $\alpha<\omega_1$ there are no infinite projective $\fin^\alpha$-mad families.
\end{corollary}
\begin{corollary}
For each $\alpha<\omega_1$ there are no infinite $\fin^\alpha$-mad families in Solovay's model.
\end{corollary}

% Bibliography
\bibliographystyle{amsplain}
\bibliography{ramsey-mad}

\end{document}